\documentclass[a4paper,10pt]{amsart}

\usepackage{amsmath}
\usepackage{amssymb}
\usepackage{graphicx}

\newcommand{\C}{\mathbb{C}}

\usepackage{amsthm}
\theoremstyle{definition}
\newtheorem*{six} {Hilbert 16th problem}
\newtheorem*{dfn1}{Definition     1}

\newtheorem*{exa1}{Example }

\newtheorem*{prop1}{Proposition  1}

\newtheorem*{pro3}{Proposition   3}

\newtheorem*{remark2}{Remark  2}
\newtheorem*{remark3}{Remark  3}

\newtheorem*{question1}{Question   }

\theoremstyle{plain}

\begin{document}

\title[Operator theory  and Geometry]{On Operator theoretical  interpretation  for  some classical problems in geometry and  differential equations  }
\author{Ali Taghavi}

\address{Faculty of Mathematics and Computer Science,  Damghan  University,  Damghan,  Iran.}
\email{taghavi@du.ac.ir}

\date{\today}


\begin{abstract}
A  consequence of the  Gauss Bonnet theorem  is interpreted in term of operator theory by Alain Connes in his  book, Non Commutative geometry. In this note we explain  in details  about his method. We
 also
introduce an operator theoretical nature for limit cycle theory.
\end{abstract}

\maketitle
\section*{Introduction}
 The  stability of the integral of  Gaussian curvature of  a  surface under  small deformation is  proved in the   book of  Alain Connes, see \cite[page 13]{NCG} . His approach is essentially based on the  stability of   value $\tau (E,E,E)$ under the one parameter deformation of projection $E$  in a $C^{*}$ algebra where $\tau$ is  a  3\_cyclic  cocycle. In this paper we explain in detail his proof.\\

 Aside of the  Gauss Bonnet theorem, we consider an important object in the theory of ordinary  differential equations, so called "limit cycles".  A  limit  cycle is  an isolated closed orbit for  planar vector field
 \begin{equation}\label{Y}
 x'=P(x,y),\;\;y'=Q(x,y)
 \end{equation}
 This vector field  sometimes is  denoted by $X=P\frac{\partial}{\partial x}+ Q\frac{\partial}{\partial y}$.\\
 Limit cycles  are the main objects of the second part of the  Hilbert 16th problem which asks:
  \begin{six}
  Does there exist  a uniform upper bound $H(n)$, depending only on $n$,  for the number of limit cycles  of $(\ref {Y})$, where $P(x,y)$  and $Q(x,y)$  are real polynomials of degree $n$?
  \end{six}
  The  vector  field $X=P\frac{\partial}{\partial x}+ Q\frac{\partial}{\partial y}$ gives us  a linear  operator $D$
  on the  space  $C^{\infty}(\mathbb{R}^{2})$, the  space of all  complex valued smooth functions on $\mathbb{R}^{2}$ as follows:
  \begin{equation}\label{Z}
  D(U)=PU_{x}+QU_{y}
  \end{equation}
  Note that the space of Schwartz functions $\mathcal{S}$ is invariant under the operator $D$ if $P$  and $Q$ are polynomial functions. Recall that $\mathcal{S}$ is the  the  subspace of   $C^{\infty}(\mathbb{R}^{2})$ consists of
  all function $f$ which all partial derivatives tends to zero at infinity.

  We  observe that the  codimension of the  rang of $(\ref{Z})$ is  an upper bound for the   number of limit cycles of $(\ref{Y})$. This  observation is true in both case that domain of $D$ is either $C^{\infty}(\mathbb{R}^{2})$ or
  the Schwartz functions $\mathcal{S}$.  So it is very important to know whether $D$ is  a fredholm operator or  a semi fredholm operator.  Among all differential operators,  the elliptic operators are the most well behaved operators, from  the view of fredholm theory. But $(\ref{Z})$ defines   a  first order differential operator $D$with real  coefficient which obviously is not  an elliptic operator. In fact an n\_th order differential operators on $\mathbb{R}^{k}$   with real  coefficients is not  elliptic if $n$ is an odd integer.

  So apparently the operator $(\ref{Z})$  is   useless, but we  have a remedy for this problem. We associate to $D$ an nth order partial differential operator $\widetilde{D}$ with polynomial  coefficient such that $\widetilde{D}$ and $D$ are  similar operators,  as two linear maps on the  space  $\mathcal{S}$ of  Schwartz functions. The  advantage of $\widetilde{D}$ is that it has the chance to be an elliptic operator when  $n$ is an even number while the  first order operator $D$ is never elliptic. The  similarity of $D$  and $\widetilde{D}$ obviously  implies that codimension of rang $D$ is equal to the  codimension of rang $\widetilde{D}$. On the  other  hand the  codimension of rang $D$ is  an upper bound for the  number of  limit cycles of $(\ref{Z})$.\\
   We have two reason for that we restrict $(\ref{Z})$ to Schwartz functions. The first reason is that we loos the similarity of $D$ and $\widetilde{D}$ on whole space $C^{\infty}(\mathbb{R}^{2})$. The second reason is that the Schwartz functions are very important subspace of Sobolov spaces which are the  main place to apply fredholm index theory. Fortunately, as we said above, the  Schwartz space $\mathcal{S}$ is invariant under $(\ref{Z})$ when $P$ and $Q$  are polynomials functions.
 \section*{The  Gauss Bonet theorem  and operator theory}
Let $S$  be  a compact  surface in $\mathbb{R}^{3}$. The  classical Gauss Bonnet theorem says that
\begin{equation}
\iint_S  \mathcal{K} dS=\mathcal{X}(S)
\end{equation}
where  $\mathcal{K}$ is the gaussian  curvature of $S$ and $\mathcal{X}$ is the Euler characteristic of $S$. This is an amazing identity because the left side is a geometric quantity but the right side is a topological one. The topology of surface $S$ do not change with an small deformation. So $\iint_S  \mathcal{K} dS$ is  stable under small deformation. In this section we explain in detail  the  cyclic cocycles interpretation for this stability as described in
 \cite[page 13]{NCG} \\
Let $N=(f,g,h)$ be the   Gauss normal map as a map  $N:S \rightarrow S^{2}$. Then the  following  formula can  be  found in any book on differential geometry of  surfaces:
\begin{equation}\label{D}
\iint_{S}  \mathcal{K} dS=\sum_{i} \iint_{G_{i}} Det \begin{pmatrix}
N &N_{u}  & N_{v}
\end{pmatrix}dudv=\sum_{i} \iint_{G_{i}} Det\begin{pmatrix}
f &f_{u}  &f_{v} \\
g & g_{u} &g_{v} \\
h& h_{u} &h_{v}
\end{pmatrix}
\end{equation}
where $G_{i}'s\subseteq \mathbb{R}^{2}$  are the domain of parameterizations  of $S$.\\
Note that $df=f_{u}du+f_{v}dv$,\;\; $dg=g_{u}du+g_{v}dv$  and  $dh=h_{u}du+h_{v}dv$.
By expanding the  determinant in $(\ref{D})$ in term of first column  we have
$$\iint_{S}  \mathcal{K} dS=\sum_{i} \iint_{G_{i}} fdgdh  +\sum_{i} \iint_{G{i}} gdhdf +\sum_{i} \iint_{G{i}} hdfdg$$
We replace the variables $u$ and $v$ in $\mathbb{R}^{2}$ with points of surface $S$ we obtain:
$$\iint_{S}  \mathcal{K} dS= \iint_{S} fdgdh  + \iint_{S} gdhdf + \iint_{S} hdfdg$$
Consider the differential 1\_forms $\alpha =fgdh$  and $\beta =ghdf$. Since $S$ is a two dimensional manifold without boundary then the  Stocks theorem implies that $\iint_{S} d(\alpha)=\iint_{S}d(\beta)=0$. This  shows that
\begin{equation}\label{T}
\iint_{S} fdgdh = \iint_{S} gdhdf = \iint_{S} hdfdg
\end{equation} So we have
\begin{equation}\label{curvature}
\iint_{S}  \mathcal{K} dS= 3\iint_{S} fdgdh
\end{equation}
 The  above  equation is  a motivation to define a trilinear form $\tau$ on the  algebra $C^{\infty}(S)$ with $\tau(f_{0},f_{1},f_{2})=\iint_{S} f_{0}df_{1}df_{2}$. Then $(\ref{T})$ shows that
\begin{equation}\label{cyclic}
\tau(f_{0},f_{1},f_{2})=\tau(f_{1},f_{2},f_{0})
\end{equation}
Furthermore  a simple  application of Leibnitz formula implies that $$f_{0}f_{1}d(f_{2})d(f_{3})-f_{0}d(f_{1}f_{2})d(f_{3})+ f_{0}d(f_{1})d(f_{2}f_{3})-f_{3}f_{0}d(f_{1}d(f_{2})=0$$
so we have
\begin{equation}\label{cocycle}
\tau(f_{0}f_{1},f_{2},f_{3})-\tau(f_{0},f_{1}f_{2},f_{3})+\tau(f_{0},f_{1},f_{2}f_{3})-\tau(f_{3}f_{0},f_{1},f_{2})=0
\end{equation}
The relation $(\ref{cyclic})$ is a cyclic property and the relation $(\ref{cocycle})$ is a cocycle property. That is every $\tau$ which satisfies in $(\ref{cocycle})$ is  a cocycle of a  complex with respect to the Hochschild coboundary map.
So  we  have  the  following definition
\begin{dfn1}
A trilinear form on  an algebra $A$  is called a 3\_cyclic cocycle if it satisfies in $(\ref{cyclic})$  and $(\ref{cocycle})$.
\end{dfn1}
Let $a$ be  an element of a unital algebra $A$   and $i,j\in\{1,2\}$.  By $\delta_{ij}(a)$ we mean a 2 by 2 matrices  in $M_{2}(A)$ which has only one nonzero entry $a$ at $i,th$ row and $j,th$ column. Every 2 by 2 matric can be written  as a sum of the  above special matrices. Now assume that $\tau$ is a 3\_cyclic cocyle on $A$. We extend $\tau$ on $M_{2}(A)$ by $$\tau(\delta_{ij}(a),\delta_{i'j'}(a'),\delta_{i''j''}(a''))=\tau(a,a',a'')trace \;(\delta_{ij}(1)\delta_{i'j'}(1)\delta_{i''j''}(1))$$
where "trace" is the  standard trace on $M_{2}(\mathbb{C})$.\\
Now consider the 3\_cyclic cocycle $\tau(f,g,h)=\iint_{S} fdgdh$ on the  algebra $A=C^{\infty}(S)$. A long but straightforward computation shows that the extension $\tau$ on $M_{2}(A)$ satisfies
$\tau(E,E,E)=\lambda \iint_{S} fdgdh$ where $\lambda$ is an independent constant and E has the following  form
\begin{equation}\label{E} E=1/2\begin{pmatrix}
1-h  & f+ig\\
f-ig & 1-h
\end{pmatrix}
\end{equation}
Here $f,g\; \text{and}\; h$ are in $C^{\infty}(S)$ with $f^{2}+g^{2}+h^{2}=1$.
We observe that \begin{equation}\label{projection}
E=E^{2}=E^{*}
\end{equation}
 where $*$ is the natural involution of $M_{2}(C^{\infty}(S))$. An element of a $C^{*}$ algebra which satisfies $(\ref{projection})$ is called a projection.
Now considering $(\ref{curvature})$ we have the  following formula for the integral of the  gaussian curvature on a surface $S$  with the Gauss normal map $N=(f,g,h)$:
\begin{equation}
\iint_{S} \mathcal{K}=\lambda \Large{\tau(E,E,E)}
\end{equation}
where $E$ is as in $(\ref{E})$ and $\lambda$ is  a universal constant.\\

Now assume that $S$ and $S'$ are two  surfaces with the Gauss normal maps $N=(f,g,h)$ and $N'=(f',g',h')$ and Gaussian curvatures
$\mathcal{K}$ and $\mathcal{K'}$, respectively. We  denote by $E$ and $E'$, the  corresponding projections of $S$  and $S'$ as defined in $(\ref{E})$.\\
 Let $\phi:S\rightarrow S'$ be a diffeomorphism between two surfaces.
Put $\widetilde{f}=f'\circ \phi$,\;$\widetilde{g}=g'\circ \phi$ and $\widetilde{h}=h'\circ \phi$ then
as a consequence of the integral version of the change of variable formula we have $\iint_{S} \widetilde{f}d\widetilde{g}d\widetilde{h}=\iint_{S'} f'dg'dh'=\iint_{S'}\mathcal{K'}=\tau(E',E',E')$. Suppose that $S'$ is a result of an small deformation of $S$. This implies that $S'$ is diffeomorphic to $S$ and $N'$ is an small deformation of $N$.
In particular  $E'$  is an small perturbation of $E$.
In the  situation of small perturbation, the following proposition says  that $\tau(E,E,E)=\tau(E',E',E')$. Thus
$$\iint_{S} \mathcal{K}=\iint_{S'} \mathcal{K'}$$
This equality shows the stability of the integral of curvature of surfaces under deformations.\\
Before we state the proposition, we recall the homotopy equivalent and unitary equivalent of projections in a $C^{*}$ algebra.  Two projection $E$ and $E'$ are homotopy  equivalent if  there is a continuous curve $E(t)$ of projections, $t\in [0,\;1]$, such that $E(0)=E,\;E(1)=E'$. $E$ is unitary equivalent to $E'$ if there is a unitary element $u$  such that $E'=uEu*$. If $E(t)$ is continuous curve of projections, then there is  a continuous curve of unitaries $u(t)$ with $u(0)=1$ such that
\begin{equation}\label{formula}
E=E(t)=u(t)E(0){u(t)}^{*}
\end{equation}

 see  \cite[Corollary 5.2.9]{OLSEN}
\begin{prop1}
Let $E(t)$  be a differentiable curve of projections in a $C^{*}$ algebra $A$. Suppose that $\tau$ is a 3\_cyclic cocycle on $A$. Then $\tau(E(t),E(t),E(t)$ is a constant complex number, for all t.
\end{prop1}

\begin{proof}
Let $E(t)$ be a differentiable curve of projections with $E(0)=E$. Then
\begin{eqnarray*}
\dot{E}=dE/dt=\lim_{t \to 0} \frac{E(t)-E}{t}     \\
 =\lim_{t \to 0} \frac{u(t)Eu(t)^{*}-E }{t}\\
 = \lim_{t \to 0} \frac{u(t)E-Eu(t)}{t}u^{*}(t)\\
 =\lim_{t \to 0}\frac{u(t)E-Eu(t)}{t}
 \end{eqnarray*}
This  shows that $\dot{E}$ is in the form of  a commutator. That is $\dot{E}=EX-XE$ for a continuous curve $X(t)$ in the algebra $A$. Now the time derivative of trilinear form $\tau(E,E,E)$ is equal to $\tau(\dot{E},E,E)+\tau(E,\dot{E},E)+\tau(E,E,\dot{E})=0$. Because we substitute  $\dot{E}$ by  $XE-EX$ and use $(\ref{cyclic})$ and $(\ref{cocycle})$. Thus $\tau(E,E,E)$ is a constant scalar.
\end{proof}

 \section*{The  number of limit cycles via differential operator theory}
 Let $A$ be  a subalgebra of $C^{\infty}(\mathbb{R}^{2})$ . We say that $A$ separates compact submanifolds of $\mathbb{R}^{2}$
 if for every finite number of disjoint  smooth closed curves $\gamma_{1},\ldots ,\gamma_{n}$ there is a real valued $f\in A$
 such that $f(x_{i}) \neq f(x_{j})$ for all $x_{i} \in \gamma_{i}$ and all $x_{j}\in \gamma_{j}$, $i\neq j$. Obviously
  $C^{\infty}(\mathbb{R}^{2})$ or algebra $\mathcal{S}$ of  Schwartz function can separates smooth closed curves in the plane. Let $X$ be a smooth vector field as in $(\ref{Y})$, in the introduction. This vector field defines the first order differential operator $D$, as in $(\ref{Z})$. Then we have:\\
  \begin{pro3}
  Let $A$ be  a self adjoint  subalgebra of $C^{\infty}(\mathbb{R}^{2})$  which separates compact submanifolds of $\mathbb{R}^{2}$.
  Assume that $A$ is invariant under the  the differential operator $D$ in $(\ref{Z})$. Then the number of closed orbits of $(\ref{Y})$ is less than or equal to the  codimension of the rang of $D$ restricted to $A$.
  \end{pro3}
  \begin{proof}
  Let $\gamma(t)=(x(t),y(t))$ be  a periodic solution for $(\ref{Y})$. Then for every real valued $g\in \C^{\infty}(\mathbb{R}^{2})$
  we have $D(g)(\gamma(t))=\LARGE{\frac{d(g\circ \gamma)}{dt}(\gamma(t))}$.  Since $g\circ \gamma$ is a periodic function, its time derivative must be vanished at some time $t_{0} \in \mathbb{R}$. So $D(g)$ must vanished on at least one point of each periodic orbit $\gamma$. Let $(\ref{Y})$ has $n$ periodic solutions(closed orbits)  $\gamma_{1},\ldots ,\gamma_{n}$.
  Assume that $f\in A$ be a real valued function which separates $\gamma_{i}$'s. Suppose that for  real numbers $\lambda_{0},\ldots,\lambda_{n-1}$, $\sum_{i=0}^{n-1} \lambda_{i}f^{i}$ is in the rang of $D$. We  show that $\lambda_{i}=0$ for $i=0,1,\ldots,n-1$. Consider the polynomial $p(z)=\sum_{i=0}^{n-1} \lambda_{i}z^{i}$. If  $p$ is a nonzero polynomial,  then there exist an $i\in \{1,2,\ldots,n\}$  such that $p(z)\neq 0$ for all $z\in \gamma_{i}$, since $\gamma_{i}$'s are n disjoint closed curves. This shows that there is no real valued function $g\in A$ such that$D(g)=p(f)$. since $A$ is  self adjoint  and the  differential operator $D$  has real coefficients, we conclude that every nontrivial combination $\sum_{i=0}^{n-1} \lambda_{i}f^{i}$  with complex coefficients is not in the rang of $D$, since f is  a real valued function. This shows that the  dimension of $A/ \text{rang of D}$,  as  a complex vector space, is at least
  $n$
  \end{proof}
  \begin{remark2}
 We observe that the codimension of the rang of the differential operator $D$ play an important role in counting the number of closed orbits of a vector field. The  codimension of the range of an operator is closely related to the concept of fredholm index. Recall that the fredholm index of an operator $D$ is equal to $\dim \ker D  -\dim Coker D$
where $\dim coker D$ is the  codimension of the rang of $D$. On the other hand, the most relevant differential operators which can be  a fredholm operator are elliptic operators. An elliptic operator is a a differential operator for which its principle part satisfies in certain nondependence condition, for  definition of elliptic operators see \cite{RENARDY}. Now the problem is that a differential operator $D$ of first order (of odd order) can not be an elliptic operator, if the  coefficients of $D$ are real valued functions.\\
 Now assume that the $P(x,y)$  and $Q(x,y)$ in $(\ref{Z})$ are polynomials of degree $n$. Then the space of Schwartz functions $\mathcal{S}$ is invariant under $D$ in $(\ref{Z})$.  Let $\mathcal{F}$ be the Fourier transform
which is a bijection from $\mathcal{S}$ to $\mathcal{S}$. Then obviously two operators $D$ and
$ \widetilde{D}=    \mathcal{F}^{-1}D\mathcal{F}$ are similar operators, hence their range have equal codimensions.
Moreover $\widetilde{D}$ is an n\_th order differential operator with polynomial coefficient. There are examples which shows that   $\widetilde{D}$ can be an elliptic operator.

  \end{remark2}

  \begin{exa1}
  Let $P(x,y)$ and $Q(x,y)$ in $(\ref{Y})$ be two polynomials of degree two which last homogeneous part are
  $P_{2}(x,y)=ax^{2}+bxy+cy^{2}$  and $Q_{2}(x,y)=dx^{2}+exy+fy^{2}$. Then using the standard  formula for the Fourier transform, which can be found for example in \cite[page 154]{RENARDY}, we obtain that the principle part of the differential operator
  $\mathcal{F}^{-1}D\mathcal{F}$ is equal to
  $$ (ax+dy)U_{xx}+(bx+ey)U_{xy}+(cx+fy)U_{yy}$$
  Now by a very simple computation we observe that this second order operator is an elliptic operator for some  appropriate coefficient of $P_{2}$  and $Q_{2}$.
  \end{exa1}
  \begin{remark3}
  As we see in the above example, in some special case ,we  can obtain an elliptic operator from a polynomial vector field on the plane. But in reality, the elliptic operators are fredholm operators only on compact manifolds, while $\mathbb{R}^{2}$ is not compact. On the other hand,  there is an standard method, so called the Poincare compactification which carries a polynomial vector fields to an analytic vector field on $S^{2}$. For  more information on poincare compactification see  \cite{LDP}. Now it would be interesting to carry the corresponding elliptic operator $\widetilde{D}=\mathcal{F}^{-1}D\mathcal{F}$ from $\mathbb{R}^{2}$ to $S^{2}$.
  What smooth vector bundle on $S^{2}$ is the best framework for consideration of such (possible) elliptic operator?

  \end{remark3}
  
  We end the paper with the following question:
  \begin{question1}
  The linear operator $(\ref{Z})$ can be considered as a pure algebraic linear map on $\mathbb{R}[x,\;y]$, the polynomial ring in two variables, when $P$ and $Q$ are polynomial. Is there an example of such operator such that the codimension of the range of $D$ is a finite number different from 0 and 1. Note that there are trivial example of codimension equal to 0, 1 and $\infty$. These codimension occur for the following vector field, respectively:
  $\frac{\partial}{\partial x}$, $x \frac{\partial}{\partial x} +y\frac {\partial}{\partial y}$ and $x\frac{\partial}{\partial x}$.
  
  \end{question1}

\end{document}